\documentclass[a4paper]{amsart}

\usepackage{amsmath,amssymb}

\newtheorem{theorem}{Theorem}
\newtheorem{corollary}[theorem]{Corollary}
\theoremstyle{definition}
\newtheorem*{definition*}{Definition}
\theoremstyle{remark}
\newtheorem*{remark*}{Remark}

\newcommand{\Prob}{{\bf Pr}}
\newcommand{\Diag}{\mbox{\rm Diag}}
\newcommand{\Ins}{\mbox{\rm Ins}}
\newcommand{\Perm}{\mbox{\rm Perm}}
\newcommand{\nat}{\mathbb N}
\newcommand{\rat}{\mathbb Q}

\begin{document}

\title{Cantorian Tableaux and Permanents}
 \author[S.Brlek]{Sre\v{c}ko Brlek}
 \thanks{LaCIM, Universit\'e du Qu\'ebec \`a Montr\'eal,
   Montr{\'e}al (QC) CANADA H3C 3P8, with the support of NSERC (Canada), 
   {\tt brlek@lacim.uqam.ca}}
 \author[M.Mend\`es France]{Michel Mend\`es France}
 \thanks{Corresponding author. D\'epartement de 
   math\'ematiques, Universit\'e Bordeaux I, 
   {\tt mmf@math.u-bordeaux.fr}}
 \author[J.M.Robson]{John Michael Robson}
 \thanks{LaBRI, Universit\'e Bordeaux I, 
   {\tt mike.robson@labri.fr}}
 \author[M.Rubey]{Martin Rubey}
 \thanks{LaBRI, Universit\'e Bordeaux I, Research financed by
   EC's IHRP Programme, within the Research Training Network \lq\lq Algebraic
   Combinatorics in Europe\rq\rq, grant HPRN-CT-2001-00272. 
   {\tt martin.rubey@labri.fr}}

\begin{abstract}
  This article could be called ``theme and variations'' on Cantor's celebrated
  diagonal argument. Given a square $n \times n$ tableau $T=\left(
    a_i^j\right)$ on a finite alphabet $A$, let $L $ be the set of its
  row-words.  The permanent $\Perm(T)$ is the set of words
  $a_{\pi(1)}^1a_{\pi(2)}^2\cdots a_{\pi(n)}^n$, where $\pi$ runs through
  the set of permutations of $n$ elements. Cantorian tableaux are those for
  which $\Perm(T)\cap L=\emptyset.$ Let $s=s(n)$ be the cardinality of $A$. We
  show in particular that for large $n$, if $s(n) <(1-\epsilon) n/\log n$ then
  most of the tableaux are non-Cantorian, whereas if $s(n) >(1+\epsilon) n/\log
  n$ then most of the tableaux are Cantorian. We conclude our article by the
  study of infinite tableaux. Consider for example the infinite tableaux whose
  rows are the binary expansions of the real algebraic numbers in the unit
  interval. We show that the permanent of this tableau contains exactly the set
  of binary expansions of \underline{all} the transcendental numbers in the
  unit interval.
\end{abstract}
\maketitle
\section{Definitions}\label{sec:definitions}
Let $A=\{a_1,a_2,\dots,a_s\}$, $s\ge 2$ be a finite alphabet and let $T$ be a
square $n\times n$ tableau
\begin{equation*}
T=\begin{pmatrix}
  a_1^1 & a_1^2 & . & . & . & a_1^n\\
  a_2^1 & a_2^2 & . & . & . & a_2^n\\
  . & . & & & & .\\
  . & . & & & & .\\
  a_n^1 & a_n^2 & . & . & . & a_n^n
  \end{pmatrix}, 
  a_i^j \in A.  
\end{equation*}
Each row $l_i=a_i^1 a_i^2 \cdots a_i^n$ is considered as a word of length $n$.
The sequence of rows is denoted by $\overline{L}$ and the set of distinct
row-words is denoted by $L$. It contains at most $n$ words.

The permanent of an $n\times n$ matrix $(a_i^j)$ defined on a ring is
\begin{equation*}
  \sum_{\pi \in S_n} a_{\pi(1)}^1a_{\pi(2)}^2\cdots a_{\pi(n)}^n,
\end{equation*}
where the summation is over the set of permutations of the $n$ elements. Very
naturally we define the permanent of the tableau $T$ to be the set of words
\begin{equation*}
  \Perm(T) = \bigcup_{\pi \in S_n} a^1_{\pi(1)}a^2_{\pi(2)}\cdots a^n_{\pi(n)}.
\end{equation*}
This set contains in particular the diagonal word
\begin{equation*}
  \Diag(T)=a_1^1 a_2^2\cdots a_n^n.
\end{equation*}
Note that the permanents of two tableaux that differ only in the order of their
rows are the same.  It may be useful to note that
\begin{equation*}
  \Perm(T)=\{\Diag(T')|\text{$T'$ a tableau obtained from $T$ by permuting its
             rows}\}.
\end{equation*}
Cantor's famous diagonal argument is based on the comparison of the set of rows
of an infinite tableau with its diagonal \cite{cantor}.  Here, at least in the
beginning of the present article, we are mostly concerned with finite $n\times
n$ tableaux and their diagonals. The last section is dedicated to infinite
tableaux.
\begin{definition*}
  A tableau is {\it Cantorian} if none of its row-words appear in $\Perm(T)$.
  In symbols
  \begin{equation*}
    L \cap \Perm(T) = \emptyset.   
  \end{equation*}
\end{definition*}
Here are some examples on the two letter alphabet $A=\{a,b\}$:
\begin{equation*}
  \begin{pmatrix}a&b\\b&a\end{pmatrix},\quad
  \begin{pmatrix}
    a&a&b&a&a&b\\
    b&b&a&b&b&a\\
    a&b&a&b&a&b\\
    b&a&b&a&b&a\\
    b&b&b&a&b&b\\
    a&a&a&b&a&a
  \end{pmatrix},\quad 
  \begin{pmatrix}
    a&b&a\\b&a&b\\b&b&b
  \end{pmatrix}.
\end{equation*}
The first one is clearly Cantorian.  Verifying that the second one also is
Cantorian seems to be a formidable task, since $\Perm(T)$ consists of $6!=720$
words. In Section~\ref{sec:sufficient} we present a simple condition which
establishes that the tableau is Cantorian. On the other hand, the third one is
not since $bbb \in L \cap \Perm(T)$.
\begin{remark*}
Note that the property of being Cantorian is invariant under permutation of
rows and columns, and, given any bijection on the alphabet, replacing all
entries of a column by their image under this bijection. To illustrate the
latter, consider the following two tableaux:
\begin{equation*}
  \begin{pmatrix}
    a&a&b&b&c\\
    a&a&b&b&c\\
    a&a&b&b&c\\
    b&b&a&a&d\\
    b&b&a&a&d
  \end{pmatrix}
  \text{ and }
  \begin{pmatrix}
    a&a&a&a&a\\
    a&a&a&a&a\\
    a&a&a&a&a\\
    b&b&b&b&b\\
    b&b&b&b&b  
  \end{pmatrix}.
\end{equation*}
While it might be difficult to see whether the first of them is Cantorian or
not, it is clear that the second is in fact Cantorian. However, it differs from
the first one only by exchanging $a$'s and $b$'s in column three and four, and
writing $a$ instead of $c$ and $b$ instead of $d$ in column five.  Hence, both
of the two tableaux must be Cantorian.
\end{remark*}

How can we calculate a permanent? There is actually an induction formula which
is similar to the one for determinants. Given an $n\times n$ tableau $T$, let
$T_i^j$ be the $(n-1)\times(n-1)$ tableau obtained by deleting row $i$ and
column $j$.  Let
\begin{equation*}
  \Ins^{\rightarrow j}\big(a,\Perm(T_i^j)\big)
\end{equation*}
denote the set of words obtained by inserting the letter $a$ at the $j$-th
place of each word in $\Perm(T_i^j)$.
\begin{theorem}\label{induction}
  For all $i \in \{1,2,\dots,n\}$
  \begin{equation*}
    \Perm(T)=\bigcup_{j=1}^n \Ins^{\rightarrow j}\big(a_i^j,\Perm(T_i^j)\big).  
  \end{equation*}
\end{theorem}
\begin{proof}
  The proof is obvious.
\end{proof}
\begin{corollary}\label{sparse}
  Let $T$ be an $n\times n$ tableau over the alphabet $A$.  Suppose a letter --
  say $a$ -- occurs $n^2-n+1$ times or more often in $T$.  Then $T$ is
  non-Cantorian. More specifically
  \begin{equation*}
    a^n \in L \cap \Perm(T). 
  \end{equation*}
  If $a$ occurs only $n^2-n$ times, the result need not be true.
\end{corollary}
\begin{proof}
  If $T$ contains no letter other than $a$ the result is trivial.  If not, we
  argue by induction on $n$. If $n=1$ the result is trivially true.  Otherwise
  let $T$ be an $n\times n$ tableau which contains at least $n^2-n+1$
  occurrences of the letter $a$. There is at least one row, say the $i$-th,
  which contains no letter other than $a$.  Since $T$ contains at least one
  letter different from $a$, there is a column $j$ with at most $n-1$
  occurrences of $a$. Thus $T_i^j$ contains at least
  $n^2-n+1-n-(n-2)=(n-1)^2-(n-1)+1$ occurrences of $a$.  By hypothesis
  \begin{equation*}
    a^{n-1} \in \Perm(T_i^j)    
  \end{equation*}
  and therefore
  \begin{equation*}
    a^n \in  \Ins^{\rightarrow j}\big(a,\Perm(T_i^j)\big).    
  \end{equation*}
  By Theorem~\ref{induction} we have that $a^n \in \Perm(T)$.
  
  On the other hand, it is easy to see that the following $n\times n$ tableau
  with $n^2-n$ occurrences of $a$ is Cantorian:
  \begin{equation}\label{example}
    \begin{pmatrix}
      a&a&.&.&.&a\\
      a&a&.&.&.&a\\
      .&.& & & &.\\
      .&.& & & &.\\
      .&.& & & &.\\
      a&a&.&.&.&a\\
      b&b&.&.&.&b
    \end{pmatrix}.
  \end{equation}
\end{proof}

The problem remains to characterize $n\times n$ Cantorian tableaux. We shall
not be able to give a definite answer to this question. We shall however
provide a sufficient condition which implies that our first two examples are
indeed Cantorian.

\section{A sufficient condition}\label{sec:sufficient}
Let $\Sigma$ be the family of all maps $\sigma:A\rightarrow A$ with no fixed
points: $\forall a\in A: \sigma(a)\not=a$.  If $aa'a''\cdots$ is a finite or
infinite word on $A$, we define
\begin{equation*}
  \sigma(aa'a''\cdots) = \sigma(a)\sigma(a')\sigma(a'')\cdots.
\end{equation*}
Let $\overline{\sigma}=(\sigma_1,\sigma_2,\dots,\sigma_n)\in\Sigma^n$. Recall
that $\overline{L}=(l_1,l_2,\dots,l_n)$ is the sequence of row-words of the
tableau $T$. Define
\begin{equation*}
  \overline{\sigma}\overline{L} = 
  (\sigma_1l_1,\sigma_2l_2,\dots,\sigma_nl_n)
\end{equation*}
and let $\overline{\sigma} L$ denote the set of all distinct $\sigma_il_i$.
Finally, let $\overline{\sigma}T$ be the $n\times n$ tableau whose row-words
are $\sigma_1 l_1,\sigma_2 l_2,\dots,\sigma_n l_n$, that is
\begin{equation*}
  \overline{\sigma}T=
  \begin{pmatrix}
    \sigma_1a_1^1 & \sigma_1a_1^2 & . & . & . &\sigma_1a_1^n\\
    . & . & & & & .\\
    . & . & & & & .\\ 
    . & . & & & & .\\
    \sigma_na_n^1 & \sigma_na_n^2 & . & . & . &\sigma_na_n^n
  \end{pmatrix}.
\end{equation*}
\begin{theorem}\label{sufficient}
  Let $\overline{\sigma}\in \Sigma^n$. Then
  \begin{equation}\label{Un}
    \Perm(T)\cap\overline{\sigma} L= \emptyset 
  \end{equation}
  and
  \begin{equation}\label{Deux}
    \Perm(\overline{\sigma}T)\cap L= \emptyset.
  \end{equation}
\end{theorem}
\begin{proof} 
  \eqref{Un} Suppose $\Perm(T)\cap \overline{\sigma} L\not= \emptyset$. Thus
  there is a permutation $\pi\in S_n$ and an index $i\in \{1,2,\dots,n\}$ such
  that
  \begin{equation*}
    a_{\pi(1)}^1a_{\pi(2)}^2\cdots a_{\pi(n)}^n =
    \sigma_i(l_i)=\sigma_i(a_i^1)\sigma_i(a_i^2)\cdots \sigma_i(a_i^n).  
  \end{equation*}
  Let $j$ be $\pi^{-1}i$. Comparing the $j$-th letter on both sides we obtain
  $a_i^j=\sigma_i(a_i^j)$, which contradicts $\sigma_i\in\Sigma$.  \eqref{Deux}
  Suppose now $\Perm(\overline{\sigma}T)\cap L \not= \emptyset$.  There is a
  permutation $\pi$ and an index $i$ such that
  \begin{equation*}
    \sigma_{\pi(1)}(a_{\pi(1)}^1)\sigma_{\pi(2)}(a_{\pi(2)}^2)
    \cdots\sigma_{\pi(n)}(a_{\pi(n)}^n) = a_i^1 a_i^2\cdots a_i^n.  
  \end{equation*}
  As above let $j=\pi^{-1}(i)$ and consider the $j$-th letter on both sides. We
  obtain
  \begin{equation*}
    \sigma_i(a_i^j) = a_i^j,
  \end{equation*}
  which again contradicts $\sigma_i\in\Sigma$.
\end{proof}

\begin{corollary}\label{cor4}
  If $L = \overline{\sigma} L$, then the tableau is Cantorian.
\end{corollary}

It is easy to see that Theorem~\ref{sufficient} holds in fact -- mutatis
mutandi -- for fixed-point free relations $\sigma$, i.e. relations with
$a\not\in\sigma(a)$ for $a\in A$. As a corollary we obtain

\begin{corollary}\label{mike}
  If for every row $i$ there exists a row $i'$ such that $a_i^j\not= a_{i'}^j$
  for all $j$, then the tableau is Cantorian.
\end{corollary}

Note that the only admissible map $\sigma\in\Sigma$ on a two letter alphabet
$\{a,b\}$ is $\sigma(a)=b,\sigma(b)=a$. Thus, in this case the
Corollaries~\ref{cor4} and \ref{mike} state the same thing.

All our previous examples of Cantorian tableaux are of the type described by
the corollary above. There are, however, other Cantorian tableaux, as the
following examples show:
\begin{equation*}
  \begin{pmatrix}
    a&a&a&a\\
    a&a&a&a\\
    b&b&b&a\\
    b&b&b&b
  \end{pmatrix},\text{ or }
  \begin{pmatrix}
    a&a&a&a&a\\
    a&a&a&a&a\\
    b&b&b&a&b\\
    b&b&b&b&a\\
    b&b&b&b&b
  \end{pmatrix}.
\end{equation*}

The following fixed point theorem is another consequence of
Theorem~\ref{sufficient}.  Before stating it we need to extend the concept of a
permanent.  Let $W=\{w_1,w_2,\dots,w_m\}$ be a set of $m$ distinct words of
length $n\ge m$.  If $m<n$, repeating some of these words, we can obtain an
$n\times n$ tableau $T$. Ignoring permutations, there are actually
$\binom{n-1}{m-1}$ ways to construct such a tableau containing all the words of
$W$.  To each one of these $T$ corresponds a permanent $\Perm(T)$.  We define
$\Perm(W)$ as the union of all the $\binom{n-1}{m-1}$ permanents.

\begin{corollary}
  Let $W=\{w_1,w_2,\dots,w_m\}$ be a set of $m\le n$ words, each of length $n$.
  Suppose that $\sigma : A\rightarrow A$ is a map such that $\sigma(W) \subset
  W$ and $W\cap \Perm(W)\not= \emptyset$.  Then there exists a letter $a\in A$
  such that $\sigma(a)=a$.
\end{corollary}
\begin{proof} Obvious from Theorem~\ref{sufficient} with
  $\overline{\sigma}=(\sigma,\sigma,\dots,\sigma)$.
\end{proof}

\section{Counting Cantorian tableaux}

Let us denote by $c(n,p)$ the number of Cantorian tableaux of size $n\times n$
over the alphabet $\{a,b\}$ having exactly $p$ occurrences of $b$.  Clearly
$c(n,p)$ has a symmetric distribution with respect to $p$, that is to say,
$c(n,p)= c(n, n^2-p)$. We also have the following computational evidence:
\begin{table}[h]
  \small
  \begin{equation*}
    \begin{array}{r|r|r|r|r|r|r|r|r|r|r|r|r|r|r|}
n\backslash p&1&2&3&4&5&6&7&8&9&10&11&12&13&14\\ \hline
2    &0     & 4     & 0  & 0 & & &  && & & & && \\
3    &0     & 0     & 3     & 9     & 9   & 3  & 0  & 0 & 0  & & & & &   \\
4    &0     & 0     & 0     & 4     & 0     & 112     & 384     & 744     & 384& 112  & 0     & 4  &0 &\ldots        \\
5    &0     & 0     & 0     & 0     & 5     & 0     & 0     & 275     & 1650     & 5960     & 14250     & 22100     & 22100 &\ldots   \\
    \end{array}
  \end{equation*}
  \caption{$c(n,p)$}
\end{table}

These numbers suggest the following result:
\begin{theorem} 
  Let $c(n,p)$ be the number of Cantorian tableaux over the alphabet $\{a,b\}$
  with exactly $p$ occurrences of the letter $b$. We have
\[   c(n,p) = 
   \begin{cases}
     0 &\text{ for } p<n,\\
     n &\text{ for } p=n\geq 3,\\
     0 &\text{ for } p=n+1\text{ and } n\geq 4,\\
     0 &\text{ for } p=n+2\text{ and } n\geq 5.
   \end{cases}
\]
\end{theorem}
\begin{remark*}
  For $p=n+3$ and $n\geq3$, the following tableau is Cantorian:
  \begin{equation*}
    \begin{pmatrix}
      \\
      \\
      \\
       &     & &b&b&b\\
      b&\dots&b&b&b&b 
    \end{pmatrix},
  \end{equation*}
  where all the entries which are not indicated are $a$'s. Hence, $c(n,n+3)$
  does not vanish.
\end{remark*}
\begin{proof}
  The $n$ Cantorian tableaux for $p=n$ are those obtained by permuting the rows
  of the tableau displayed in \eqref{example} in the proof of
  Corollary~\ref{sparse}. We will show that there are no others by considering
  several cases.
\begin{itemize}
\item If $p<n$, Corollary~\ref{sparse} applies and we are done. 
\item If $n\leq p\leq n+2$ and if there is a row $b^n$, direct inspection of the
  possible cases shows that the tableau is Cantorian if and only if $p=n$.
  Up to permutation of rows and columns, the possible cases are
  \begin{gather*}
    \begin{pmatrix}
      \\ \\ \\
      b&\dots&b&b 
    \end{pmatrix},
    \begin{pmatrix}
      \\ \\
       &     & &b\\
      b&\dots&b&b 
    \end{pmatrix},\\
    \begin{pmatrix}
      \\ \\
       &     &b&b\\
      b&\dots&b&b 
    \end{pmatrix},
    \begin{pmatrix}
      \\
       &     & &b\\
       &     & &b\\
      b&\dots&b&b 
    \end{pmatrix},
    \begin{pmatrix}
      \\
       &     &b&\\
       &     & &b\\
      b&\dots&b&b 
    \end{pmatrix},
  \end{gather*}
  where all the entries which are not indicated are $a$'s.
\item If $n\leq p\leq n+2$ and there is no row $b^n$, but a row $a^n$ we claim
  that $a^n$ is in the permanent of $T$: let $i$ be the row with the greatest
  number of $b$'s. Because there is a row without any $b$'s and the total
  number of $b$'s is at least $n$, there are at least two letters $b$ in this
  row. Since there is no row $b^n$, this row contains at least one letter $a$,
  say in column $j$.
  
  Now we consider the tableau $T_i^j$.  Clearly, $a^{n-1}$ is a row of $T_i^j$.
  Furthermore, note that $b^{n-1}$ cannot be a row of $T_i^j$: this would be
  possible only if there were two rows $b^{n-1}$ in $T$.  Since the number of
  $b$'s in $T$ is $p$, we would have $p\geq 2(n-1)$, and hence $n\leq p-n+2$,
  which contradicts the bounds we assumed for $n$.
  
  We proceed by induction on $n$: since $T_i^j$ contains at most $p-2$ letters
  $b$, we have by hypothesis that $a^{n-1}$ is in the permanent of $T_i^j$.
  Since $\Ins^{\rightarrow j}\big(a,\Perm(T_i^j)\big)$ is a subset of the
  permanent of $T$, we are done.
\item If $n\leq p\leq n+2$ and there is neither a row $b^n$, nor a row $a^n$,
  we have to consider two subcases:
  \begin{itemize}
  \item There is a column $j$ with at least two $b$'s, with one of them being
    the only one in its row. Let $i$ be one of the other rows having a $b$ in
    column $j$. Consider the tableau $T_i^j$. Clearly, it contains a row
    $a^{n-1}$. Because $T$ has no row $a^n$, the reduced tableau $T_i^j$ cannot
    contain a row $b^{n-1}$: otherwise we had $p\geq n-1+2+n-3=2n-2$ which
    contradicts the bounds we assumed for $n$.
    
    Now the statement established in the previous case applies to $T_i^j$, and
    by Theorem~\ref{induction} we obtain that row $i$ is in the permanent of
    $T$, which contains $\Ins^{\rightarrow j}\big(b,\Perm(T_i^j)\big)$.
  
  \item Otherwise, by permuting rows and columns, the tableau can be
    represented as
    \begin{equation}
      \begin{pmatrix}
        b \\
        &\ddots\\
        &      &b\\
        &      & &b&b\\
        &      & &b&b
      \end{pmatrix},    
    \end{equation}
    where all the entries which are not indicated are $a$'s. Clearly, this
    tableau is non-Cantorian.
  \end{itemize}
\end{itemize}
\end{proof}

\section{An Algorithm for Enumerating Cantorian Tableaux}\label{sec:algorithm}
Let $C(n,s)$ be the number of $n\times n$ Cantorian tableaux on an s letter
alphabet.  Computing $C(n,s)$ is obviously quite cumbersome even for $s=2$.  As
a first improvement over simple-minded calculation of the permanent followed by
checking whether the intersection with the set of row-words is nonempty, we
have the following:

\begin{theorem}\label{polynomial}
  It is possible to test whether a given tableau is Cantorian or not in
  polynomial time.
\end{theorem}
\begin{proof}
  To test whether a given $n\times n$ tableau $T$ over any alphabet is
  Cantorian or not, we proceed as follows: for each row $k$ we transform the
  tableau into a bipartite graph $G_k$, with each of the two parts having $n$
  vertices. The $n$ \lq top\rq\ vertices $\{x_1,x_2,\dots,x_n\}$ correspond to
  the rows of the tableau, the $n$ \lq bottom\rq\ vertices
  $\{y_1,y_2,\dots,y_n\}$ correspond to the columns of $T$.  There is an edge
  connecting the \lq top\rq\ vertex $x_i$ with the \lq bottom\rq\ vertex $y_j$
  if and only if the entries in column $j$ in row $i$ and in row $k$ are the
  same, i.e., if $a_i^j=a_k^j$.
  
  If $G_k$ contains a perfect matching, then the tableau $T$ cannot be
  Cantorian.  Otherwise, we proceed with the next line of the matrix. If there
  is no perfect matching for any of the rows of the tableau, it must be
  Cantorian.
  
  Since it is possible to find a perfect matching in a given bipartite graph in
  time $O(n^{2.5}/\sqrt{\log n})$, where $n$ is the number of vertices of the
  graph -- see \cite[Theorem~1]{AltBlumMehlhornPaul} --, we see that this
  procedure is polynomial in time.
\end{proof}

For very small $n$, we can inspect all $n\times n$ tableaux and check whether
they are Cantorian.  Since, even for a two letter alphabet, the number of
tableaux is $2^{n^2}$, this rapidly becomes infeasible as a method for
determining their number.  We have computed the number of Cantorian tableaux of
sizes up to $9$ over the two letter alphabet $\{0,1\}$ as follows:

First we consider only tableaux whose last row contains only ones.  Any tableau
can be converted into one of this form by zero or more ``column flips''
consisting of changing every element in a chosen column. We obtain the total
simply by multiplying the number obtained by $2^n$.

For row $k$ consider the bipartite graph $G_k$ as defined in the proof of
Theorem~\ref{polynomial}. Since this graph always contains the edge
$(x_k,y_n)$, if the graph $G'_k$ obtained by omitting $x_k$ and $y_n$ has a
perfect matching, then so does $G_k$. Thus, by generating all directed
bipartite graphs $G'_n$ with each part consisting of $n-1$ vertices and without
a perfect matching, we can find all Cantorian tableaux of size $n$ by $n$.  For
$n=9$ the number of such graphs is ``only'' $1256511813403160577$. We improve
on this idea as follows:
\begin{itemize}
\item {\bf Canonicity.} We define an equivalence relation on the graphs and
  generate only one canonical instance from each equivalence class. We choose
  the equivalence relation so that the size of the equivalence class is easy to
  determine and each graph in an equivalence class produces the same number of
  Cantorian tableaux.  The number of canonical graphs with no perfect matching
  for $n=9$ is only $11446766661$.
\item{\bf Last column conditions.} Given a graph $G'_n$ without a perfect
  matching, that is, given the first $n-1$ columns of a tableau, we can quickly
  determine the number of ways to fill the last column. The condition that the
  complete tableau is Cantorian is a conjunction of conditions of the form
  $a_i^n\not=a_{i'}^n$, which are easy to obtain. The number of ways of
  choosing the final column so as to give a Cantorian tableau is then either
  zero, if the conditions are inconsistent, or $2^{n-i}$, where $i$ is the
  number of independent conditions. This calculation is significantly faster
  than generating and testing all $2^{n-1}$ final columns with $a_n^n=1$.
\item{\bf Skeletons.} We generate the tableaux corresponding to the graphs
  $G'_n$ as follows: first we fix only a \emph{skeleton} -- a small portion of
  the entries of the tableau, about a third in each row and column. Then we
  decide upon the value of the remaining entries a row at a time. For all rows
  $k$ which are already completed, we test whether $G'_k$ has a perfect
  matching. If this is the case we can discard the generated tableau.  In
  practise, this happens at a very early stage, when the value of only few
  entries of the tableau has been fixed.
\end{itemize}

With all these improvements, the calculation for $n=9$ takes a little less than
an hour. The results can be found in Table~2 below.
\section{Asymptotics}\label{sec:asymptotics}
\begin{table}[h]
  \footnotesize
  \begin{equation*}
    \begin{array}{r|r|r|r|r|r}
      n\backslash s  &                         2&              3&
            4&            5&             6\\ \hline
                    2&                  1\cdot 2^2&     4\cdot 3^2&
   9\cdot 4^2&  16\cdot 5^2&   25\cdot 6^2\\
                     &            2.5\cdot 10^{-1}&\sim 4.44\cdot 10^{-1}&
5.625\cdot 10^{-1}&6.4\cdot 10^{-1}&\sim 6.94\cdot 10^{-1}\\ \hline  
                    3&                  3\cdot 2^3&   188\cdot 3^3&
1863\cdot 4^3&9264\cdot 5^3&32075\cdot 6^3\\
                     &      \sim 4.69\cdot 10^{-2}&\sim 2.58\cdot 10^{-1}&
\sim 4.55\cdot 10^{-1}&\sim 5.93\cdot 10^{-1}&\sim 6.87\cdot 10^{-1}\\ \hline  
                    4&                109\cdot 2^4&100144\cdot 3^4&&&\\
                     &      \sim 2.66\cdot 10^{-2}&\sim 1.88\cdot 10^{-1}&&&\\ \hline
                    5&               2765\cdot 2^5&&&\\
                     &      \sim 2.64\cdot 10^{-3}&&&\\ \hline
                    6&             324781\cdot 2^6&&&\\
                     &      \sim 3.02\cdot 10^{-4}&&&\\ \hline
                    7&           37304106\cdot 2^7&&&\\
                     &      \sim 8.48\cdot 10^{-6}&&&\\ \hline
                    8&        13896810621\cdot 2^8&&&\\
                     &      \sim 1.93\cdot 10^{-7}&&&\\ \hline
                    9&      5438767247337\cdot 2^9&&&\\
                     &      \sim 1.15\cdot 10^{-9}&&&\\ \hline
                   10&6889643951630251\cdot 2^{10}&&&\\
                     &     \sim 5.56\cdot 10^{-12}&&&
      \end{array}
  \end{equation*}
  \caption{Number and proportion of Cantorian tableaux of size $n$ on 
    alphabets of size $s$.}
\end{table}

Consider the $16$ tableaux of size $2\times 2$ over the alphabet $\{a,b\}$.
Direct inspection shows that among them, only $4$ are Cantorian. There are
$2^{n^2}$ tableaux of size $n\times n$.  It is reasonable to guess that among
them there is only a small proportion of Cantorian tableaux.  Let $C(n)=C(n,2)$
be the number of Cantorian tableaux and $N(n)$ the number of non-Cantorian
tableaux over $\{a,b\}$. We have the following explicit bounds:
\begin{equation*}
  N(n)\ge 2^{n^2-n+1}\text{ and }
  C(n) >
  \begin{cases}
    2^{\frac{1}{2}n^2}&\text{if $n$ is even,}\\
    2^{\frac{1}{2}n(n-1)}&\text{if $n$ is odd.}
  \end{cases}
\end{equation*}
In particular, the lower bound for $C(n)$ is obtained as follows.  Suppose $n$
is even. Choose arbitrarily the entries in the first $n/2$ rows of an $n\times
n$ tableau $T$. Complete $T$ by adjoining $n/2$ rows obtained by interchanging
the letters $a$ and $b$ in the first $n/2$ rows of $T$. By Corollary~\ref{cor4}
this tableau is necessarily Cantorian, and there are $2^{\frac{1}{2}n^2}$ such
tableaux.  Obvious modifications establish the case where $n$ is odd.

Therefore
\begin{equation*}
\frac{1}{2}\leq \liminf_{n\to\infty} \frac{\log C(n)}{\log 2^{n^2}}\leq \limsup_{n\to\infty}
\frac{\log C(n)}{\log 2^{n^2}}\leq 1.  
\end{equation*}
The question arises whether the limit actually exists, and if so, what is its
value? Here are a few values of the ratio
\begin{equation*}
  \frac{\log C(n)}{\log 2^{n^2}} = .5,.509,.673,.657,.675, .656,.651,.632,
  .626 
\end{equation*}
for respectively $n= 2,3,4,\ldots, 10$, computed from Table~2, which displays
the values of the number $C(n,s)$ of $n\times n$ Cantorian tableaux on an
alphabet of size $s$, for some values of $n$ and $s$.

The proportion of Cantorian tableaux on a fixed $s$-letter alphabet tends to
$0$ as the size of the tableaux increases, as is suggested by Table~2 and
established by the following theorem:
\begin{theorem}\label{asymptotic}
  Let $C(n,s)$ be the number of $n\times n$ Cantorian tableaux on an alphabet
  of size $s=s(n)$. If $s< n/(\log n+\log\log n+r_n)$ where $r_n$ is any
  sequence which grows without bound, then
  \begin{equation*}
    \lim_{n\rightarrow \infty}C(n,s)/s^{n^2} = 0. 
  \end{equation*}
  If on the other hand, $s> n/(\log n-\log\log n-\epsilon)$ for any $\epsilon
  >0$, then
  \begin{equation*}
    \lim_{n\rightarrow \infty}C(n,s)/s^{n^2} = 1.    
  \end{equation*}
\end{theorem}
\begin{proof}
  A tableau $T$ is certainly non-Cantorian if $l_n\in \Perm(T)$. This in turn
  is certainly the case if there is a permutation $\pi\in S_{n-1}$ such that
  $a_i^{\pi(i)}=a_n^{\pi(i)}$ for all $i<n$. If $\pi$ consists of a single
  cycle, then the following directed graph $G$ has a (directed) Hamiltonian
  cycle: $G=(V,E)$, where $V=\{1,2,\dots,n-1\}$ and $(i,j)\in E$ if and only if
  $a_i^j=a_n^j$. Note that $G$ is derived from the graph $G^\prime_n$ in
  Section~\ref{sec:algorithm} by directing all its edges from \lq top\rq\ to
  \lq bottom\rq\ and then identifying vertices $x_i$ and $y_i$ for
  $i\in\{1,2,\dots,n-1\}$.
  
  This graph $G$ for a randomly chosen $n\times n$ tableau is simply
  $D_{n-1,1/s}$, a random directed graph on $n-1$ vertices where each possible
  edge has probability $1/s$ of occurring. The probability that such a graph is
  Hamiltonian is well known to tend to $1$ as $n$ tends to infinity
  \cite{frieze}, as long as the alphabet size $s(n)$ grows with $n$ but is
  bounded by $s(n)<n/(\log n +\log \log n+r_n)$ where $r_n$ is any sequence
  which grows without bound.
  
  On the other hand, if $s(n)>n/(\log n -\log \log n-\epsilon)$ for
  $\epsilon>0$, then Corollary~\ref{mike} shows that the probability that a
  random tableau is Cantorian tends to $1$. Indeed, for any two rows $i$ and
  $j$,
  \begin{align*}
    \Prob [\forall k:a_i^k \not= a_j^k] 
    &> (1-(\log n-\log\log n -\epsilon)/n)^n \\
    &> \log n(1+\epsilon')/n~ \quad\quad (\epsilon'>0)
  \end{align*}
  for $n$ sufficiently large. Hence for a given $i$,
  \begin{align*}
    \Prob [\exists j \forall k:a_i^k \not= a_j^k] 
    &> 1-(1-\log n(1+\epsilon')/n))^{n-1} \\
    &> 1-n^{-(1+\epsilon'')}~ \quad\quad(\epsilon''>0)
  \end{align*}
  for $n$ sufficiently large. We deduce that the expected number of $i$ such
  that $\nexists j \forall k:a_i^k \not= a_j^k$ is less than
  $n^{-\epsilon''}$.  Therefore the probability that there exists such an $i$
  is less than $n^{-\epsilon''}$ and thus tends to $0$.
\end{proof}

\section{Infinite Tableaux}
The definitions and results of the preceding sections extend naturally to
infinite tableaux $T=(a_i^j)$ with $i,j\in \nat$.  In particular the permanent
of $T$ is the set of infinite sequences
\begin{equation*}
  \Perm(T)=\bigcup_{\pi\in S_{\nat}}
  a_{\pi(1)}^1a_{\pi(2)}^2a_{\pi(3)}^3\cdots,
\end{equation*}
where $S_{\nat}$ is the family of all bijections $\pi: \nat\rightarrow \nat$.
In general, $\Perm(T)$ is an uncountable set.

Consider an infinite tableau $T$ over the alphabet $A$.  The $i$-th row ($i\in
\nat$) is
\begin{equation*}
  l_i= a_i^1 a_i^2 a_i^3 \cdots {}\in A^{\nat}
\end{equation*}
and the set of rows is denoted $L$ as in the finite case.  If
$\overline{\sigma}=(\sigma_1,\sigma_2,\sigma_3,\dots)\in \Sigma^{\nat}$, where
$\Sigma$ is defined as in Section~\ref{sec:sufficient}, then
Theorem~\ref{sufficient} extended to the infinite case asserts that
\begin{equation*}
  \Perm(T)\cap \overline\sigma L= \emptyset
  \text{ and }
  \Perm(\overline{\sigma}T)\cap L = \emptyset.
\end{equation*}
If moreover $L=\overline{\sigma}L$, then $T$ is Cantorian.

\begin{theorem}\label{transcendental}
  Let $L$ be a countable subset of the unit interval such that $L\supseteq
  {\rat\cap[0,1]}$ and $L+\rat=L \bmod 1$. Let $T$ be the infinite tableau
  whose rows are the expansions in base $s\ge 2$ of the numbers in $L$. Here we
  require that rational numbers $r/s^q$ for $r,q\in\nat$ should appear twice in
  $T$, once with a tail of $0$'s and once with a tail of $s-1$'s. Then $T$ is
  Cantorian: $\Perm(T)$ only contains numbers in $[0,1] \setminus L$.  The
  $s$-expansion of each number in $\Perm(T)$ contains each of the digits
  $0,1,\dots,s-1$ infinitely often. None of the digits occur periodically.
\end{theorem}
\begin{remark*}
  In fact, for $s\ge 3$ the theorem holds also if the tableau contains only one
  of the two possible expansions of rational numbers $r/s^q$ for $r,q\in\nat$,
  i.e., the expansion having an infinite tail of $0$'s.
  
  For $s=2$ the intersection $\Perm(T)\cap L$ contains nothing but the numbers
  whose expansion has an infinite tail of $0$'s. Since $0$ is such a number,
  the statement that both digits $0$ and $1$ occur infinitely often, is false
  in this setting.
\end{remark*}
\begin{proof}
  Let $l_i=a_i^1a_i^2a_i^3\cdots$ represent both the $i$-th row of $T$ and the
  $i$-th element of $L$. Assume that $s\geq 3$. In this case, we define 
  \begin{equation*}
    \pi:l_i\mapsto l_i+1/(s-1) \bmod 1.
  \end{equation*}
  Clearly, $\pi$ is a reordering of the rows. Moreover, since
  $1/(s-1)=.1111\dots$, it can be shown that the $j$-th digit of $\pi(l_i)$ is
  always different from the $j$-th digit of $l_i$. For $s=2$ we define
  $\pi:l_i\mapsto 1-l_i$, which implies the same fact. By Corollary~\ref{mike}
  we conclude that the tableau is Cantorian, i.e., $\Perm(T)$ contains only
  expansions of numbers in $[0,1] \setminus L$.
  
  Next we prove that each of the $s$-digits $0,1,\dots,s-1$ occurs infinitely
  often in every element of $\Perm(T)$.  Suppose the digit $a$ occurs only
  finitely many times in
  \begin{equation*}
    p=a_{\pi(1)}^1 a_{\pi(2)}^2\cdots.
  \end{equation*}
  Define the map $\sigma:\{0,1,\dots,s-1\}\rightarrow\{0,1,\dots,s-1\}$ by
  $\sigma(b)=a$ for all $b\not=a$ and $\sigma(a)=c$, where $c$ is any letter
  different from $a$.  Then $\sigma\in\Sigma$ and therefore $\Perm(\sigma T)$
  contains no numbers from $L$. Thus
  \begin{equation*}
    \sigma(p)=\sigma(a_{\pi(1)}^1)\sigma(a_{\pi(2)}^2)\cdots
  \end{equation*}
  is the $s$-expansion of a number in $[0,1] \setminus L$.  But from some point
  on, $p$ contains no digit $a$, so that $\sigma(p)$ has an infinite tail of
  $a$'s. This is absurd since then $\sigma(p)\in\rat$.  The same map $\sigma$
  shows that no $s$-digit of $p$ can occur periodically.
\end{proof}

\begin{remark*}
  Theorem~\ref{asymptotic} asserts that given an alphabet $A$, with a fixed
  number $s$ of elements, the probability that an $n\times n$ tableau is
  Cantorian tends to $0$ as $n$ increases. Therefore
  Theorem~\ref{transcendental} should come as a surprise. Paradoxically, the
  infinite tableau $T$ described in Theorem~\ref{transcendental} is Cantorian.
  This of course is due to the fact that the set $L$ is closed under the
  addition of rational numbers, a rather stringent condition indeed!
\end{remark*}

\begin{theorem} \label{exact}
  If $T$ is as in the statement of Theorem~\ref{transcendental} with $s=2$,
  then
  \begin{equation*}
    \Perm(T)=[0,1]\setminus L.   
  \end{equation*}
\end{theorem}
\begin{proof}
  By Theorem~\ref{transcendental} the tableau $T$ is Cantorian, thus we have
  that $\Perm(T)\subseteq [0,1]\setminus L$. The equality is established as
  follows: let $x$ be a number in $[0,1]\setminus L$.  We show how to construct
  a permutation $\pi\in S_{\nat}$ such that
  \begin{equation*}
    x=0.a_{\pi(1)}^1a_{\pi(2)}^2\cdots.
  \end{equation*}
  Writing $x^j$ for the $j$-th digit of $x$ we define
  \begin{equation*}
    \pi(j)=\min\{i|~a_i^j=x^j \text{ and } \forall j'<j: i\not=\pi(j')\}.
  \end{equation*}
  
  First we show that $\pi(j)$ is well defined: since $L$ contains $\rat$ there
  is an infinity of rows $i$ with $a_i^j=x^j$. On the other hand, there can be
  only a finite number of rows $i$ with $i=\pi(j')$ for some $j'<j$. Thus,
  the set of which we take the minimum is indeed nonempty.
  
  It is clear from the definition that $\pi$ is one-to-one. Hence it remains to
  show that for all rows $i$ there is a column $j$ such that $\pi(j)=i$.  Note
  that there is an infinity of columns $j$ with $a_i^j=x^j$, because otherwise
  $x+l_i$ would be rational, where $l_i=a_i^1 a_i^2\cdots$. This in turn cannot
  be the case, since $l_i$ is in $L$ and $x$ is not.
  
  Suppose now that for all $j$ we have $\pi(j)\not=i$.  It follows that
  $\pi(j)\leq i$ for all columns $j$ with $a_i^j=x^j$.  Since $\pi$ is
  one-to-one this can be true for only a finite number of columns, thus
  contradicting the assumption.
\end{proof}

By choosing the set of algebraic numbers in the unit interval for $L$, we
obtain the following corollary:
\begin{corollary}\label{cortranscendental}
  If the rows of $T$ consist of all the algebraic numbers in the unit interval
  represented in base 2, then $\Perm (T)$ is exactly the set of all
  transcendental numbers in the unit interval.
\end{corollary}

The condition $s=2$ was necessary in Theorem~\ref{exact} and Corollary
\ref{cortranscendental}. Indeed, for $s\ge 3$ there exist transcendental numbers with
missing digits such as the Liouville numbers
\begin{equation*}
  \sum_{n\ge 0} s^{-n!}.
\end{equation*}
Therefore, Theorem~\ref{transcendental} shows that for $s\ge 3$, $\Perm(T)$
cannot contain all transcendental numbers.

It is however true for all $s\ge 2$ that $\Perm(T)$ contains uncountably many
transcendental numbers. Indeed, suppose $l_i$ is the list of all the algebraic
numbers as above and $t_i$ is a list of some countable set of transcendental
numbers.  We show how to construct a permutation $\pi$ such that
\begin{equation*}
  a_{\pi(1)}^1 a_{\pi(2)}^2\cdots
\end{equation*}
is different from every $t_i$: note that any $l_i$ and $t_j$ differ from each
other in infinitely many positions. Let $i_1$ be the first position where $l_1$
differs from $t_1$. We choose $\pi(1)=i_1$ and we set $\pi(j)=j-1$ for all
$1<j<i_1$. Clearly $a_{\pi(1)}^{i_1}\not=t_1^{i_1}$.  Now we proceed
iteratively: when $i_k$ is known, we choose $i_{k+1}$ as the first position
after $i_k$ where $l_{i_k+1}$ differs from $t_{k+1}$. Then we define
$\pi(i_k+1)=i_{k+1}$ and $\pi(j)=j-1$ for all $j$ such that $i_k+1 < j <
i_{k+1}$. $\pi$ permutes each set $\{i_k+1,i_k+2,\dots,i_{k+1}\}$, so it is
indeed in $S_{\nat}$ and $a_{\pi(i_k+1)}^{i_{k+1}}\not=t_{k+1}^{i_{k+1}}$.

We can in fact show the stronger result that the set of real numbers not
contained in $\Perm(T)$ has measure $0$.  This is a consequence of the
following theorem.
\begin{theorem}
  Let $T$ be an infinite tableau containing the $s$-expansions of a countable
  dense subset $L$ of the unit interval. Then the measure of $\Perm(T)$ is $1$.
\end{theorem}
\begin{proof}
  We will show that the probability that $x\not\in \Perm(T)$ is less than
  $\epsilon$ for any $\epsilon>0$.  Recall that the probability that an $n$
  node random directed graph $D_{n,1/s}$ (where each possible edge has
  probability $1/s$ of being present) is Hamiltonian tends to $1$ as $n$ tends
  to infinity.  Define $n_i$ for $i>0$ to be the first $n$ such that this
  probability is greater than $1-\epsilon/2^i$ and let $N_0=0$ and
  $N_i=\sum_{j=1}^i n_j$.
  
  From the initial order of the rows $l_i$ we construct a new order $l'_i$ as
  follows: $l'_{N_i+1}$ is the first $l_j$ not already present in
  $l'_1,l'_2,\dots,l'_{N_i}$. For $N_i+2\le j \le N_{i+1}$, elements in
  positions $N_i+1,\dots, N_{i+1}$ of $l'_j$ are chosen randomly and
  independently and $l'_j$ is chosen as the first $l_k$ with these elements and
  not already present in $l'_1,,l'_2,\dots,l'_{j-1}$. Such an $l'_j$ exists
  since $L$ is dense.  Because $l_j$ is chosen at the latest in the $j$-th
  step, this defines a reordering of the rows.
  
  Now we consider the probability that $x$ can be obtained as the diagonal of a
  permutation of the rows $l'_i$ fixing each of the sets
  $\{N_i+1,N_i+2,\dots,N_{i+1}\}$. Let ${~}^iv$ be the vector consisting of
  digits $N_i+1,N_i+2,\dots,N_{i+1}$ of $x$. Let ${~}^iB$ be the square Boolean
  matrix whose entry ${~}^iB_{j}^k$ is true if and only if
  ${~}^iM_j^k={~}^iv^k$.  We claim that ${~}^iB$ has each entry true with
  probability $1/s$ and that all these probabilities are independent: this is
  true for the first row of ${~}^iB$ because $x$ was chosen randomly, and for
  all other rows because the corresponding elements of ${~}^iM$ were random.
  
  Now ${~}^iB$ is the adjacency matrix of a graph $D_{n_i,1/s}$ and we know
  that this graph has probability less than $\epsilon/2^i$ of not being
  Hamiltonian. If all the graphs are Hamiltonian there is a permutation of the
  rows $l'_i$ that consists of a cycle on each of the sets
  $\{N_i+1,N_i+2,\dots, N_{i+1}\}$ and produces $x$ on the diagonal. Hence
  \begin{equation*}
    \Prob[x\not\in \Perm(T)] < 
    \sum_{i=1}^{\infty}\epsilon/2^i =\epsilon. 
  \end{equation*}
\end{proof}
 
A similar approach allows us to consider tableaux whose rows consist of
infinite sequences on a finite set $\{0,1,\dots,s-1\}$. Rational numbers would
be replaced by ultimately periodic sequences and algebraic numbers would then
be replaced by $s$-automatic sequences \cite{allouche,chr, ckmr}. The results from
Theorem~\ref{transcendental} on remain valid with the obvious modifications.
 
\section{Outlook}
Following the remark after the definition of Cantorian tableaux in
Section~\ref{sec:definitions}, define an equivalence relation on the set of
$n\times n$ tableaux as follows: let $T'$ be equivalent to $T$, if it is
obtained from $T$ by a combination of permuting rows or columns or replacing
all entries of a column by their image under any bijection on the alphabet. It
might be interesting to count the number of resulting equivalence classes.

Taking into account the situation for base $2$ in Theorem~\ref{transcendental},
it might also be interesting to consider those tableaux $T$ where $\Perm(T)\cap
L$ equals a given set, or has a given cardinality.

Finally, we could have defined ``bi-Cantorian'' tableaux as those where
$\Perm(T)$ is disjoint both from the set of row-words and column-words. We
chose our initial definition guided by Cantor's work. Needless to say it might
well be interesting to extend our discussion to bi-Cantorian tableaux. For
example, an argument very similar to the one given at the beginning of
Section~\ref{sec:asymptotics} shows that there are at least
$2^{{\lfloor n/2\rfloor}^2}$ $n\times n$ bi-Cantorian tableaux over the alphabet
$\{a,b\}$.

\end{document}